\documentclass[10 pt, conference]{ieeeconf}  

\IEEEoverridecommandlockouts                              
\title{\LARGE \bf 
Sensitivity analysis of hybrid systems with state jumps \\
with application to trajectory tracking}	
\date{\today}

\author{Alessandro Saccon, Nathan van de Wouw, Henk Nijmeijer %
\thanks{A. Saccon, N. van de Wouw, and H. Nijmeijer are with the
  Department of Mechanical Engineering,
  Eindhoven University of Technology, Eindhoven, the Netherlands.
  {\tt\small \{a.saccon, n.v.d.wouw, h.nijmeijer\}@tue.nl}
  }%
\thanks{
  Research supported by the 
  European Union Seventh Framework Programme [FP7/2007-2013]
  under grant agreement no. 257462 HYCON2 Network of excellence.  
}
}

\usepackage{amsmath}
\usepackage{amssymb}
\usepackage{bm}       
\usepackage{graphicx}

\usepackage{color}

\newcommand {\eps} {\varepsilon}
\newcommand {\R} {\mathbb{R}}

\newcommand {\te}{\tau} 

\newtheorem{theorem}{Theorem}[section]

\newtheorem{proposition}[theorem]{Proposition}

\newtheorem{assumption}{Assumption}[section]

\begin{document}

\maketitle

\begin{abstract}
This paper addresses the sensitivity analysis 
for hybrid systems with discontinuous (jumping) state trajectories.
We consider state-triggered jumps in the state evolution, potentially accompanied by 
mode switching in the control vector field as well.
For a given trajectory with state jumps, we show how to construct 
an approximation of a nearby perturbed trajectory corresponding 
to a small variation of the initial condition and input.
A major complication in the construction of such an approximation
is that, in general, the jump times corresponding
to a nearby perturbed trajectory are not equal to those of the nominal one. 
The main contribution of this work is the development 
of a notion of error to clarify in which sense 
the approximate trajectory is, at each instant of time, 
a first-order approximation of the perturbed trajectory.
This notion of error naturally finds application in the (local) tracking problem 
of a time-varying reference trajectory of a hybrid system. 
To illustrate the possible use of this new error definition in 
the context of trajectory tracking,
we outline how the standard linear trajectory tracking control
for nonlinear systems --based on linear quadratic regulator (LQR) theory 
to compute the optimal feedback gain--
could be generalized for hybrid systems.
\end{abstract}

%
%
%

\section{Introduction}

Sensitivity analysis for dynamical systems allows to determine directly the change in a trajectory due to (small) changes in initial conditions and parameters and has proven beneficial in many aspects of the analysis of dynamical systems.
In this paper, we pursue such sensitivity analysis for a class of hybrid systems.
A hybrid system is a dynamic system that exhibits both continuous and discrete dynamic behaviors \cite{GoSaTe12B_HybridSystems}.
It is this inherent nature of hybrid systems that makes sensivity analysis for this class of systems harder 
than for nonlinear dynamics, where sensitivity analysis is well established (see, e.g., \cite[Chapters 3]{Kh02B_NLSystems}).
The theoretical framework provided by hybrid systems with (state-triggered) jumps 
is suitable to model those systems which, at certain instants of time, are subjected to rapid and abrupt changes.
Indeed, in the modeling of such systems, it is frequently convenient and valid 
to neglect the durations of these rapid changes and to assume that the changes can 
be represented by instantaneous state jumps.
Hybrid systems with state-triggered jumps are, for example, suitable 
to describe dynamical models in the area of robotics and rigid body mechanics with unilateral contact constraints \cite{LeinevandeWouw2008BOOK,BR99B_NonsmoothMechanics},
including the study of the dynamics and control of walking or juggling robots (see, e.g,
\cite{WeGrKo03J_HybridZeroDynPlanarBipedWalker, RoLeSe06J_BounceJuggling, MoGr09J_HybridInvManifoldSystemsWithImpulseEffects}).

The sensitivity result, presented in this paper, draws from the application
of classical sensitivity and perturbation theory of nonlinear systems 
(see, e.g., \cite[Chapters 3 and 10]{Kh02B_NLSystems})
combined with the use of the implicit function theorem
to compute an estimate of the unknown switching time 
at which the perturbed trajectory jumps.
These same mathematical tools have been combined 
in the investigation of the sensitivity about a nominal trajectory 
of piecewise-smooth nonlinear systems without state jumps in 
\cite{LeCaVr00J_BifurcNLDiscontinuousSys}, where the concept of {\em saltation matrix} was introduced.
One can also find them in \cite{HiPa00J_SensitivityHybridSystems}, where
although state jumps are not considered, it is recognized that part the analysis could 
be carried on even in the presence of state jumps.
In \cite[Section 6.4]{Mo01D_StabOptOpenLoopWalkRobots}, one finds an interesting discussion regarding 
the sensitivity of hybrid systems in the context of numerical optimal control of mechanical systems 
where state jumps are specifically taken into account.
This work refers to \cite[Section 2.2] {ScWiSc96C_ParEstimDiscDescriptorModel} 
that in turn refers to \cite{Kr85B_SeminalSensitivityAccordingToMombaur} (in German)
as the source of a key formula for defining the sensitivity of hybrid systems with state jumps.
We will re-establish  this key formula for the class of hybrid systems considered in this work. From now on, we will refer to the sensitivity over a jump 
as the {\em jump gain} associated with a discontinuous event and it is given in Equation \eqref{eq:Hstep} in Section \ref{sec:sensitivity}.
Another interesting discussion regarding sensitivity analysis for hybrid system with state jumps 
is presented in \cite{GaFeBa99J_SensitivityFunctionsHybridSystems}, for hybrid systems 
with jumps obtained by combining multiple differential algebraic equations (DAEs) with 
switching conditions and reset maps.
There, it is mentioned that the jump gain \eqref{eq:Hstep} should be credited in fact credited 
to the (seemingly forgotten) seminal work of Rozenvasser \cite[equation (11)]{Ro67J_GenSensitivityEqsDiscSystems}.
%
For completeness, we also mention that
a formula {\em related} to \eqref{eq:Hstep}
can be encountered in the context of
discontinuity induced bifurcations 
for hybrid systems. The interested reader is referred to
\cite[Section 2.5.1]{BeBuChKo08B_PiecewiseSmoothDynSystems}
and the work in \cite{No91J_NonPeriodMotionGrazingImpactOscillator},
where the concept of discontinuity map (for transversal intersections) is introduced.

We leave to historians the settling of the question on who discovered \eqref{eq:Hstep} first.
Here, we limit ourselves to mentioning that the jump gain \eqref{eq:Hstep} is indeed a key result that 
we, as other researcher before us, rediscovered autonomously. 
Our goal in this paper is to rigorously define what 
a first-order approximation about a nominal trajectory of a hybrid system with state jumps is on the basis of the result on the jump gain.
By doing so, a novel notion of error, which allows to locally compare a nominal trajectory with a perturbed trajectory, 
emerges. This error notion, in turn, leads naturally to obtain a trajectory tracking controller to locally stabilize 
time-varying trajectories for hybrid systems with state jumps. This is, in our view, the main novelty 
and contribution of this paper. Since  the proper understanding of the jump gain \eqref{eq:Hstep} is a key ingredient in such developments, we also include a (re)derivation of this result in this work.

Tracking control for hybrid system with state-triggered jumps
is a recent and active field of research.  
Few results exist to design a controller to make 
a hybrid system with state-triggered jumps 
track a given, time-varying, reference trajectory.
Recent techniques addressing this control problem 
both from a theoretical and an experimental viewpoint 
are provided in \cite{
MeTo01J_AsympTrackPeriodicTrajsNonsmothImpacts,
PaYu01C_ExpStudyPlanarImpactRobotManipulator,
Pa01J_CtrlContactProblemConstrainedEL,
LeWo08J_UnifConvofMonotoneMDI,
FoTeZa13J_FollowTheBouncingBall},
and
\cite{BiWoHeNi13J_TrackCtrlHSStateTriggeredJumps}.
Our interest lies on the situation, 
commonly encounter in practice, where the jump times 
of plant and reference trajectories cannot be assumed to coincide.

Aiming at developing an effective trajectory tracking controllers 
for hybrid systems with state-triggered jumps, 
we propose to investigate the effects 
on a nominal trajectory 
of variations of the control input and initial conditions.
In particular, we detail how to construct a linear approximation 
of the hybrid system about a nominal trajectory with jumps 
and then show how to use this approximation to construct 
a local trajectory tracking controller. 
To the best of our knowledge, 
the notion of linear approximation that we introduce in this paper
has not be presented before.

The non trivial aspect of the problem is how to construct
the approximation of the perturbed trajectory 
as the sum of the nominal trajectory 
and a linear term 
when the perturbed and nominal trajectories jump
at different, although close, time instants.
This difference in the jump times poses also the problem
of defining a proper notion of tracking error.

Mimicking what is done for dynamical systems
with no jumps, indeed, the most intuitive definition 
of tracking error is the difference between the nominal 
and actual state.
However, as illustrated in, e.g., the introduction of \cite{BiWoHeNi13J_TrackCtrlHSStateTriggeredJumps},
this definition has the drawback of 
exhibiting an unstable behavior in the sense of Lyapunov
in spite of the converge of the perturbed trajectory to the nominal one (away from the jump times).

For this reason, different approaches have been proposed 
in recent years to redefine the notion of tracking error for hybrid systems. 
In \cite{GaMePo08C_TrajTrackLinearHS}, e.g.,
the tracking problem has been defined in order to neglect in
the analysis the times belonging to infinitesimal intervals about the jumping times.
In \cite{BiWoHeNi13J_TrackCtrlHSStateTriggeredJumps}, a novel conceptual definition of the notion of distance between two jumping trajectories has been proposed. Moreover,
for a subclass class of hybrid system state-triggered jumps 
corresponding to mechanical systems with fully elastic impacts,
the tracking error distance has been defined as the minimum 
between the distance of state from the nominal trajectory 
and
the distance of state a mirrored version of nominal trajectory.
In \cite{BiWoHeSaNi12C_TrackingDissipativeImpacts}, a modification 
of this mirroring approach has been consider to deal with dissipative impacts.
The use of a mirror reference trajectory has been proposed 
recently also in \cite{FoTeZa13J_FollowTheBouncingBall}
for a tracking problem in polyhedral billiards: the controller 
in this case may decide to track either the real reference 
or the mirrored reference, mirrored through the billiard boundary.

We claim that the approximation proposed in this paper 
allows for a reinterpretation of the mirroring technique 
in terms of what we will call {\em extended ante- and post-event trajectories}.
Furthermore, we care to emphasize the fact the concept of  
extended ante- and post-event trajectories
allows to cope with the problem of trajectory tracking
for a hybrid-system with state-triggered jumps
where the (continuous-time) dynamics before and after the jump event 
are qualitatively different, where simply mirroring the reference 
trajectory does not appear to be the best choice.

This paper is organized as follows.
In Section \ref{sec:sensitivity}, we discuss the jump gain 
associated with a nominal trajectory of a hybrid system and 
introduce the notion of error between the nominal and perturbed trajectories 
of hybrid system with state jumps.
This notion of error is used in Section \ref{sec:trajTracking} 
to propose a linear feedback control law for local trajectory tracking 
of time-varying reference trajectories with jumps.
Conclusions are finally drawn in Section \ref{sec:conclusions}.

\section{Sensitivity analysis for hybrid systems}
\label{sec:sensitivity}

In this section, we propose a framework for the sensitivity analysis
of a jumping solution of a hybrid system with state-triggered jumps.
To focus on the complexity of the effect of such a jump on 
the sensitivity, we limit ourselves to sensitivity analysis
about one such a jump and leave the treatise
of sensitivity analysis of solutions with multiple jumps for future work.

Consider a (sufficiently) smooth time-varying control vector field 
\begin{align}
  \dot x(t) & = f^a(x, u, t), 
  \label{eq:ante_vectorfield}
\end{align}
with state $x \in \R^n$ and input $u \in \R^m$. For reasons that will 
appear clear shortly, we will refer to $f^a$ as the {\em ante-event} control vector field.
For a given initial condition $x_0$ at time $t_0$ and 
a integrable signal $\mu(t) \in \R^m$, $t \in [t_0,t_1]$,
we denote with $\alpha^a(t)$, $t \in [t_0, \te]$, 
the solution of \eqref{eq:ante_vectorfield} with input
\begin{align}
  u(t) & = \mu(t),  \quad t  \in [t_0,t_1],
  \label{eq:nominal_input}
\end{align}
up to the occurrence of a triggering {\em event} at time $\te \in [t_0,t_1]$
defined by the satisfaction of the implicit condition
\begin{align}
  g( \alpha^a (\te), \te ) = 0 ,
  \label{eq:nominal_eventcondition}
\end{align}
where $g:\R^n \times \R \rightarrow \R$ is a {\em smooth} real-valued function
We assume that, for all $t$, the level set $g(\cdot,t) = 0$
is a $n-1$ dimensional smooth manifold embedded in $\R^n$,
(a sufficient condition being $\partial g(x,t)/\partial x \neq 0$ for all $x$
such that $g(x,t) = 0$). At the event time $\te$, 
the state exhibits a {\em jump} according to 
a smooth {\em impulse} map $\Delta:\R^n \times \R \rightarrow \R^n$
\begin{align}
  \alpha^p(\te) = \alpha^a(\te) + \Delta(\alpha^a(\te),t),
  \label{eq:nominal_resetmap}
\end{align}
and subsequently it evolves according to the following {\em post-event} 
vector field 
\begin{align}
  \dot x(t) & = f^p(x, u, t) 
  \label{eq:post_vectorfield}
\end{align}
with initial condition at $x(\te) = \alpha^p(\te)$
and input \eqref{eq:nominal_input}. We indicate with 
$\alpha^p(t), t \in [\te, t_1]$ the post-event trajectory.

\begin{figure*}
\centering
\includegraphics[width=0.8\textwidth]{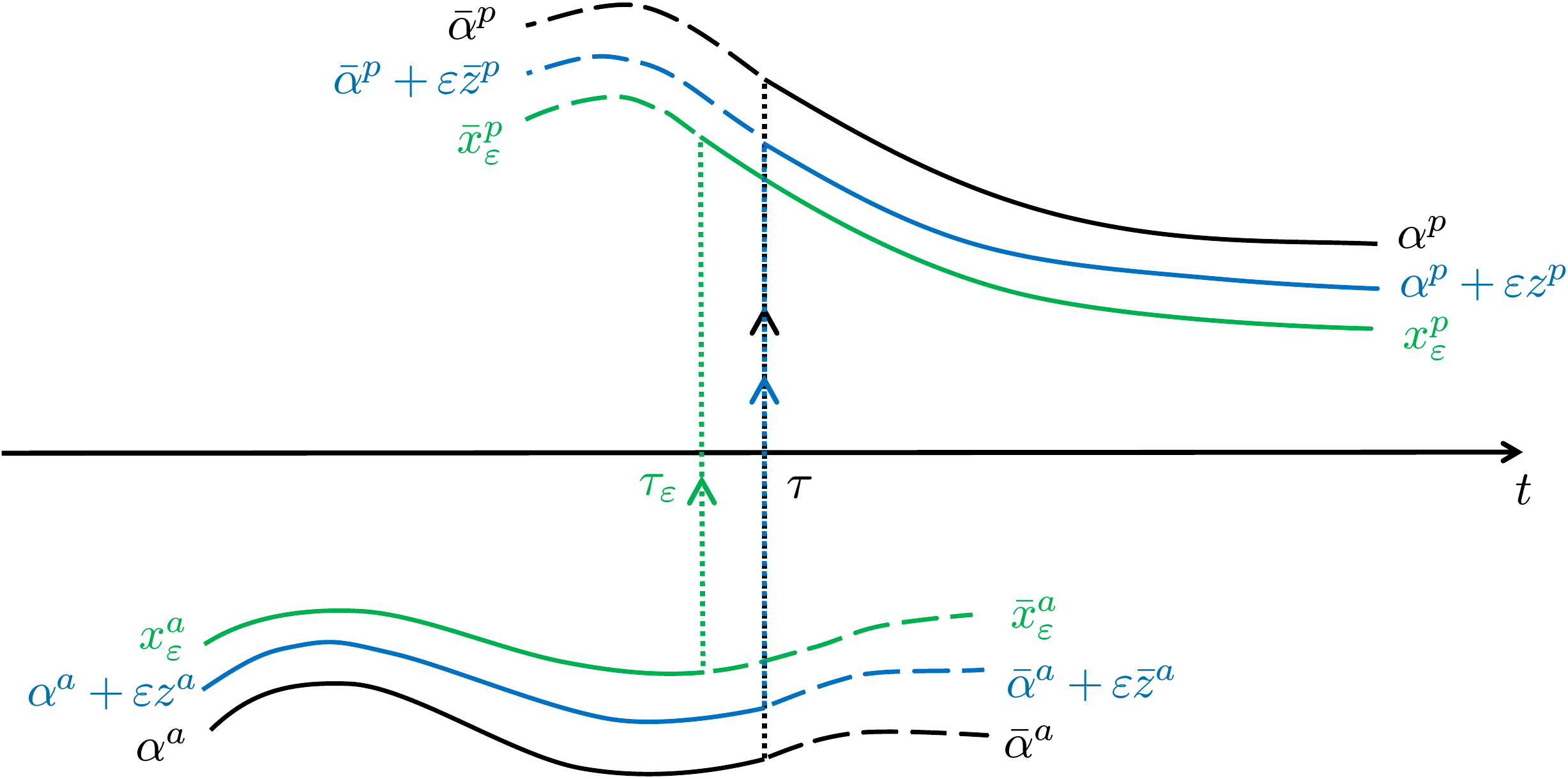}
\caption{
The ante- and post-event nominal trajectories $\alpha^a$ and $\alpha^p$,
perturbed trajectories $x_\eps^a$ and $x_\eps^p$, and the 
linear approximations $\alpha^a + \eps z^a$ and $\alpha^p + \eps z^p$.
The extended trajectories, indicated with a bar sign, are depicted using dashed lines.
}
\label{fig:alpha_plus_eps_z}
\end{figure*}

Ante- and post-event trajectories can be glued together 
by ``appending'' to ante-event trajectory $\alpha^a$
the post-trajectory $\alpha^p$, obtaining the trajectory 
\begin{align}
  \alpha(t) & := 
    \left\{ 
    \begin{array}{ll}
      \alpha^a(t), & t \in [t_0, \te) \\ 
      \alpha^p(t), & t \in [\te, t_1] .
    \end{array} 
    \right. 
    \label{eq:nominal_statetraj}
\end{align}
The state trajectory $\alpha$ is (generally) not continuous, due to the state jump caused by the impulse map $\Delta$.
Without loss of generality, $\alpha$ is by construction right continuous. 
Together with the nominal input $\mu(t)$, the state trajectory $\alpha(t)$
forms what we term the {\em nominal state-input trajectory} $\xi(t) = (\alpha(t), \mu(t))$, $t \in [t_0, t_1]$.

We are interested in defining and computing the sensitivity of $\xi(t)$ 
to small variations of the initial condition $x_0$ and input $\mu(t)$, $t \in [t_0, t_1]$. 
To this end, we will perturb the initial condition in the direction $z_0 \in \mathbb{R}^n$ 
and the nominal input curve $\mu$ in the direction  $v$, 
where $v$ denotes an integrable curve $v(t) \in \mathbb{R}^m$, $t \in [t_0, t_1]$.
Hence, the perturbed initial condition and input curve are defined, respectively, as 
\begin{align}\label{Eq:iiperturbations}
  x_\eps(t_0) & = x_0 + \eps z_0 \\
  u_\eps(t)   & = \mu(t) + \eps v(t), \quad t \in [t_0, t_1] ,  
\end{align}
with $\eps \in \R$ typically small.
The corresponding state-input trajectory will be denoted
$\xi_\eps(t) = (x_\eps(t),u_\eps(t))$, $t \in [t_0,t_1]$.
We aim at defining a notion of sensitivity 
to predict, for small values of $\eps$, 
the effect of perturbations to the initial condition and the nominal input
on a jumping solution of the hybrid system.

The perturbed trajectory $\xi_\eps = (x_\eps, u_\eps)$ is defined, similarly to the nominal trajectory $\xi = (\alpha,\mu)$,
by appending to the {\em perturbed} ante-event trajectory $x^a_\eps$ the 
the {\em perturbed} post-event trajectory $x^p_\eps$.
One cannot expect, however, that the event time $\te$ will remain constant as $\eps$ is varied.
Indeed, the perturbed event time, that we denote $\te_\eps$, is implicitly defined by the condition
\begin{align}
  g( x^a_\eps (\te_\eps), \te_\eps ) = 0 .
  \label{eq:perturbed_eventcondition}
\end{align}
The above equation is satisfied, for $\eps = 0$, by the nominal ante-event trajectory $\alpha^a$
at the nominal event time $\te$, see \eqref{eq:nominal_eventcondition}. 

In order to ensure that $\eps \mapsto \te_\eps$ is, about $\eps = 0$, a continuously differentiable function
the following assumption needs to be fulfilled by the nominal state-input trajectory $\xi$.
\begin{assumption}
The nominal state-input trajectory satisfies the following {\em transversality} condition
\begin{align}\label{Eq:transversality}
  D_1 g \cdot f^a + D_2 g \cdot 1 \neq 0 , 
\end{align}
where $g$ is evaluated at $(\alpha^a(\te),\te)$ and the ante-event
vector field $f^a$ at $(\alpha^a(\te), u(\te), \te)$.
\label{a:transversality}
\end{assumption}
The above transversality condition is a common requirement when dealing with hybrid system
with state-triggered jumps (see, e.g., \cite{MoGr09J_HybridInvManifoldSystemsWithImpulseEffects}).
Its role in the sensitivity analysis will become clear in the proof of the proposition 
presented later on in this section.

For $t = \te_\eps$, the value of $x^p_\eps$ 
is computed using the impulse map $\Delta$, similarly as done in \eqref{eq:nominal_resetmap}.
Therefore, we can formulate the following conditions that the perturbed ante- and post-event
trajectories have to satisfy:
\begin{align}
  x^a_\eps   & = x_0 + \eps z_0 ,                     &  t & = t_0 
  \label{eq:perturbed_init_cond} \\
  \dot x^a_\eps     & = f^a(x^a_\eps, u_\eps, t),   &  t & \in [t_0, \te_\eps]
  \label{eq:perturbed_ante_event} \\
  x^p_\eps   & = x^a_\eps + \Delta(x^a_\eps,t),       &  t & =    \te_\eps   
  \label{eq:perturbed_at_event}  \\
  \dot x^p_\eps     & = f^p(x^p_\eps, u_\eps, t),   &  t & \in [\te_\eps, t_1]. 
  \label{eq:perturbed_post_event}
\end{align}
One important observation is that, although naturally defined 
in the time intervals $[t_0, \te_\eps]$ and $[\te_\eps, t_1]$, 
both the ante- and post-event trajectories $x^a_\eps$ and $x^p_\eps$ can individually
be extended over the whole time interval $[t_0, t_1]$
by, respectively, forward and backward time integration starting from $\te_\eps$.   
We will denote those extensions as ${\bar x}^a_\eps(t)$ and ${\bar x}^p_\eps(t)$, 
$t \in [t_0,t_1]$. 
These extended ante- and post-event perturbed trajectories satisfy
\begin{align}
  {\bar x}^a_\eps             & = x_0 + \eps z_0  & t & = t_0 , 
  \label{eq:init_cond} \\
  \dot {\bar x}^a_\eps        & = f^a({\bar x}^a_\eps, u_\eps, t),         &  t & \in [t_0, t_1],
  \label{eq:ante_event} \\
  {\bar x}^p_\eps             & = {\bar x}^a_\eps + \Delta({\bar x}^a_\eps,t),  &  t & =    \te_\eps,   
  \label{eq:at_event}  \\
  \dot {\bar x}^p_\eps        & = f^p({\bar x}^p_\eps, u_\eps, t),         &  t & \in [t_0,t_1], 
  \label{eq:post_event}
\end{align}
where the perturbed event time $\te_\eps$ is implicitly defined by the condition
\begin{align}
  g( {\bar x}_\eps^a(\te_\eps), \te_\eps ) = 0 .
  \label{eq:g_eq_zero}  
\end{align}
This apparently innocuous extension is the cornerstone 
to understand the sensitivity differential equation. 

Figure~\ref{fig:alpha_plus_eps_z} gives an indication on why in general
we cannot expect to be able to write the perturbed trajectory $x_\eps$ as
\[
  x_\eps(t) = \alpha(t) + \eps z(t) + o(\eps)
\]
with $z$ the solution to an appropriately defined time-varying 
linear system. The obstacle is represented by the difference 
in the event times for $x_\eps$ and $\alpha$. 
The following proposition shows how to overcome this difficulty, 
defining the extended ante- and post-event linearization trajectories
${\bar z}^a$ and ${\bar z}^p$ and defining the expansion about the
extended ante- and post-event trajectories ${\bar\alpha}^a$ and 
${\bar\alpha}^p$ instead of simply about $\alpha$.


\begin{proposition}
\label{prop:approximation}
Consider a nominal state-control trajectory $\xi(t) = (\alpha, \mu)(t)$, 
$t~\in~[t_0,t_1]$, of system \eqref{eq:ante_vectorfield}, \eqref{eq:nominal_eventcondition}, \eqref{eq:nominal_resetmap}, \eqref{eq:post_vectorfield}
with nominal event time $\te \in [t_0,t_1]$ and associated 
extended ante- and post-event trajectories ${\bar\alpha}^a(t)$ and ${\bar\alpha}^p(t)$, $t \in  [t_0,t_1]$. Adopt Assumption \ref{a:transversality}.
The perturbed state trajectory $x_\eps(t)$, $t \in [t_0,t_1]$, corresponding to perturbations in the initial condition and input as in \eqref{Eq:iiperturbations}, satisfies
\begin{align}  
  x_\eps(t) & = 
  \left\{ 
  \begin{array}{ll} 
    {\bar\alpha}^a(t) + \eps {\bar z}^a(t) + o(\eps), & t < \te_\eps \\
    {\bar\alpha}^p(t) + \eps {\bar z}^p(t) + o(\eps), & t \geq \te_\eps .
  \end{array}
  \right.
  \label{eq:xeps_alpha_plus_z}
\end{align}
where the extended ante- and post-event linearization trajectories ${\bar z}^a(t)$ 
and ${\bar z}^p(t)$, $t \in [t_0, t_1]$, are computed as
\begin{align}
        {\bar z}^a    & = z_0,                           & t & = t_0 \label{eq:za_z0}\\
   \dot {\bar z}^a    & = A^a(t) {\bar z}^a + B^a(t) v   & t & \in [t_0, t_1] \label{eq:za_Aaza_plus_Bav} \\
        {\bar z}^p    & = {\bar z}^a + H(\te) {\bar z}^a & t & = \te \label{eq:zp_za_plus_Hza}\\
   \dot {\bar z}^p    & = A^p(t) {\bar z}^p + B^p(t) v   & t & \in [t_0, t_1] \label{eq:zp_Apzp_plus_Bpv}
\end{align}
where
\begin{align}
  A^s(t) & = D_1 f^s({\bar\alpha}^s(t),\mu(t),t)    \\
  B^s(t) & = D_2 f^s({\bar \alpha}^s(t),\mu(t),t)  
\end{align}
with $s = \{a,p\}$ and
\begin{align}
  H(\te)      & = \frac { f^+ - f^- - \dot \Delta^- }{ \dot g^- } D_1 g^- + D_1\Delta^- ,  
  \label{eq:Hstep}
\end{align}
where
\begin{align}
 \alpha^p(\te)      & = \alpha^a(\te) + \Delta(\alpha^a(\te),\te)                    \\
 f^+                & = f^p(\alpha^p(\te), \mu(\te), \te),                
 \label{eq:fplus}                                                                    \\
 f^-                & = f^a(\alpha^a(\te), \mu(\te), \te),                           
 \label{eq:fminus}                                                                   \\
 \dot \Delta^-      & = D_1\Delta^- \cdot f^- + D_2\Delta^- \cdot 1                  
 \label{eq:dDeltaminus}                                                              \\ 
 \dot g^-           & = D_1 g^- \cdot f^- + D_2 g^- \cdot 1                      
 \label{eq:dot_g_minus}                                                              \\
 D_k \Delta^-       & = D_k \Delta(\alpha^a(\te), \te),              & k & = \{1,2\}  
 \label{eq:DDeltaminus}                                                              \\ 
 D_k g^-            & = D_k g(\alpha^a(\te), \te),                   & k & = \{1,2\}. 
 \label{eq:D_g}
\end{align}
\end{proposition}
\vspace{0.2cm}
\begin{proof}
The extended ante- and post-event perturbed trajectories satisfy
\eqref{eq:init_cond}-\eqref{eq:post_event} where 
the perturbed event time $\te_\eps$ is implicitly defined by the condition
\eqref{eq:g_eq_zero}. On the basis of Assumption \ref{a:transversality}, the implicit function theorem allows to conclude the existence of a unique 
event time $\te_\eps$ for small values of $\eps$. 
Furthermore, as \eqref{eq:g_eq_zero} is identically zero for each $\eps$ in a neighbourhood of zero,
the derivative of \eqref{eq:g_eq_zero} with respect to $\eps$ allows to obtain a linear 
approximation of $\te_\eps$ as a function of $\eps$,
namely to compute $\te_0' := \partial \te_\eps / \partial \eps |_{\eps = 0}$
and approximate $\te_\eps$ as
\begin{align}
  \te_\eps & = \te + \eps \, \te_0' + o(\eps) .
  \label{eq:teps_expansion}
\end{align}
The formula to compute $\te_0'$ will be obtained at the end of the proof.

The sensitivity equations for the extended ante- and post-event trajectories ${\bar x}^a$ 
and ${\bar x}^p$ are the standard (see, e.g., \cite[Chapter 10]{Kh02B_NLSystems}) sensitivity equations given by 
\eqref{eq:za_Aaza_plus_Bav} and \eqref{eq:zp_Apzp_plus_Bpv}.
The missing link is how to relate $\bar z_p$ to $\bar z_a$, that is, 
	to show that ${\bar z}_p$ is indeed reinitialized as in \eqref{eq:zp_za_plus_Hza} at the {\em nominal} event time $\te$ 
using the linear map $H$ given by \eqref{eq:Hstep}. 

Expanding ${\bar x}^p_\eps(\te_\eps)$ in series with respect to $\eps$ 
results in
\begin{align}
  {\bar x}^p_\eps(\te_\eps) 
& = 
  {\bar \alpha}^p(\te_\eps) + \eps {\bar z}^p(\te_\eps) + o(\eps) \notag \\ 
& = 
  {\alpha}^p(\te) + \eps {\dot {\alpha}}^p(\te) \te_0' + \eps {\bar z}^p(\te_\eps) + o(\eps) \notag \\
& = 
  {\alpha}^p(\te) + \eps ( {\dot \alpha}^p(\te) \te_0' + {\bar z}^p(\te) ) + o(\eps) .
  \label{eq:alpha_p_series}  
\end{align}
To obtain the above expression, ${\bar \alpha}^p(\te_\eps)$ has been approximated 
  by linear extrapolation using the value and the time derivative of ${\bar \alpha}^p$ at time $\tau$
  and $\te_\eps$ has been approximated using \eqref{eq:teps_expansion}.
  We have then discarded the terms of order higher than one in $\eps {\bar z}^p(\te_\eps)$.
A similar expansion can be computed for ${\bar x}^a_\eps(t_\eps)$ obtaining
\begin{align}
  {\bar x}^a_\eps(\te_\eps) 
& = 
  {\alpha}^a(\te) + \eps ( {\dot \alpha}^a(\te) \te_0' + {\bar z}^a(\te) ) + o(\eps) .
  \label{eq:alpha_a_series}
\end{align}  
Using \eqref{eq:alpha_a_series} and \eqref{eq:teps_expansion},
$\Delta({\bar x}^a_\eps(\te_\eps),\te_\eps)$ appearing in \eqref{eq:at_event}
can be expanded as 
\begin{align}
  \Delta({\bar x}^a_\eps(\te_\eps), \te_\eps)
& = 
  \Delta({\alpha}^a(\te),\te) \notag\\ 
  & \hspace{-5ex}+ \eps \, D_1 \Delta({\alpha}^a(\te),\te) \cdot (\dot {\alpha}^a(\te) \, \te_0' + {\bar z}^a(\te) ) \notag\\ 
  & \hspace{-5ex}+ \eps \, D_2 \Delta({\alpha}^a(\te),\te) \cdot \te_0' + o(\eps) .
  \label{eq:Delta_series}
\end{align}
Using \eqref{eq:alpha_a_series} and \eqref{eq:teps_expansion},
  \eqref{eq:g_eq_zero} can be expanded as 
\begin{align}
  g( {\bar x}^a_\eps(\te_\eps), \te_\eps )
  & = \eps \, D_1 g ( \alpha^a(\te), \te) \cdot (\dot {\alpha}^a(\te) \, \te_0' + {\bar z}^a(\te) )  \notag \\
  & + \eps \, D_2 g ( \alpha^a(\te), \te) \cdot \te_0'  + o(\eps).
\end{align}
As the above expression is identically zero for every $\eps$, 
we get 
\begin{align}
  \te_0' 
& = 
 - \frac
 { D_1 g ( \alpha^a(\te), \te) \cdot {\bar z}^a(\te) } 
 { D_1 g ( \alpha^a(\te), \te) \cdot \dot {\alpha}^a(\te) + D_2 g ( \alpha^a(\te), \te) \cdot 1 } 
 \notag \\ 
& =  
  - \frac{D_1 g^- \cdot {\bar z}^a(\te)}{{\dot g}^-}  
  \label{eq:tep}
\end{align}
where ${\dot g}^-$ and $D_1 g^-$ are defined, respectively, as in \eqref{eq:dot_g_minus} and \eqref{eq:D_g}.
The expression for \eqref{eq:tep} is valid as long as ${\dot g}^-$ is different from zero, that is,
if Assumption \ref{a:transversality} is satisfied.

Making use of the series expansions \eqref{eq:alpha_p_series}, \eqref{eq:alpha_a_series}, and \eqref{eq:Delta_series},
we can match the first-order terms of \eqref{eq:at_event}, obtaining, at the event time $\te$, 
\begin{align}
  \dot {\alpha}^p  \te_0' + {\bar z}^p
  & =
  \dot {\alpha}^a \te_0' + {\bar z}^a \notag \\
  & \hspace{-10ex} 
  +
  D_1 \Delta({\alpha}^a,\te) \cdot (\dot {\alpha}^a \, \te_0' + {\bar z}^a ) 
  + 
  \, D_2 \Delta({\alpha}^a,\te) \cdot \te_0' .  
\end{align}
The above expression can be rearranged as
\begin{align}
  {\bar z}^p
& =
  {\bar z}^a   
  - (   f^+
      - f^-
      - \dot\Delta^-
    ) \, \te_0'  
  + 
    D_1 \Delta^- \cdot {\bar z}^a
\end{align}
where $f^+$, $f^-$, $\dot\Delta^-$, and $D_1 \Delta^-$
are defined, respectively, as in \eqref{eq:fplus}, \eqref{eq:fminus}, 
\eqref{eq:dDeltaminus}, and \eqref{eq:DDeltaminus}.
Substituting in the above equation the expression for $\te_0'$ given in \eqref{eq:tep},
we obtain \eqref{eq:zp_za_plus_Hza} and in particular 
the expression for $H$ provided in \eqref{eq:Hstep}.
This concludes the proof of the proposition.
\end{proof}
\noindent
{\bf Remark.} Note that the reset map \eqref{eq:zp_za_plus_Hza} is linear in ${\bar z}^a$ 
and the reset occurs at the nominal event time $\te$.\hfill$\Box$

As mentioned in the introduction, 
\eqref{eq:Hstep} is a (uncommonly) known expression
in the context of numerical optimal control 
\cite{ScWiSc96C_ParEstimDiscDescriptorModel},
\cite{Mo01D_StabOptOpenLoopWalkRobots}, 
and 
parametric sensitivity for hybrid systems
\cite{Ro67J_GenSensitivityEqsDiscSystems},
\cite{GaFeBa99J_SensitivityFunctionsHybridSystems}.
It is also strictly related to equation (57) that appears 
in \cite{HiPa00J_SensitivityHybridSystems}
and can also be interpreted 
as a generalization for piecewise-smooth nonlinear systems with state jumps 
of the saltation matrix introduced in \cite{LeCaVr00J_BifurcNLDiscontinuousSys} 
(in \cite{LeCaVr00J_BifurcNLDiscontinuousSys}, $\Delta$ is identically equal to zero as there is no state jump).

Our contribution lies in the use of extended ante- and post-event trajectories 
to achieve the $o(\eps)$ approximation presented in \eqref{eq:xeps_alpha_plus_z}. 
To the best of our knowledge, this approximation and 
the use the extended anti- and post-event trajectories is new. 

It is worth mentioning that \cite[Appendix A]{HiPa00J_SensitivityHybridSystems}
``suggests a procedure for refining the estimate of the perturbed trajectory''
in a neighbourhood of the event time. The need for this refining is due to the difference
between the nominal and perturbed trajectory event times (in our notation, $\te$ and $\te_{\eps}$).
The approximation in \cite{HiPa00J_SensitivityHybridSystems} differs from the one that we propose 
as it is not an $o(\eps)$ approximation of the nominal trajectory (uniformly in $t$) due to the lack of use 
of the extended anti- and post-event trajectories.

The approximation \eqref{eq:xeps_alpha_plus_z} 
is key to address the problem of (local) trajectory tracking for systems with state-triggered jumps, 
as we discuss in the following section.

\section{Trajectory tracking of a time-varying reference trajectory}
\label{sec:trajTracking}

Let $\xi = (\alpha, \mu)$ be a nominal state-input trajectory 
for a hybrid system characterized by ante- and post-event vector fields $f^a$
and $f^p$. As done in the previous section, $\te$ will indicate 
the nominal event time and ${\bar \alpha}^a$ and ${\bar \alpha}^p$ 
the extended ante- and post-event trajectories, respectively.
 
Consider the following state feedback control law
\begin{align}
  u = 
  \left\{
    \begin{array}{ll}
      \mu + K(t) ( {\bar\alpha}^a - x ), & \text{before event detection}  \\
      \mu + K(t) ( {\bar\alpha}^p - x ), & \text{after event detection}
    \end{array}
  \right.  
  \label{eq:u_cl}
\end{align}
where $K(t) \in \mathbb{R}^{m \times n}$ is a time-varying matrix gain to be designed.
In \eqref{eq:u_cl}, by {\em event detection} we mean the satisfaction of the condition
$g(x(t),t)=0$ for the current value of the state at time $t$ (assumed not to be equal to
the nominal event time $\te$).
As mentioned in the example of tracking control for a bouncing ball proposed in
\cite{BiWoHeNi13J_TrackCtrlHSStateTriggeredJumps},
event detection is not strictly needed to implement a switching feedback law as \eqref{eq:u_cl}.
Indeed, due to the discontinuity in the nominal state trajectory, 
an equivalent result is obtained by simply choosing, 
between ${\bar\alpha}^a - x$ and ${\bar\alpha}^p - x$, the one with minimum norm.

\noindent{\bf Remark.} 
Strictly speaking, in \cite{BiWoHeNi13J_TrackCtrlHSStateTriggeredJumps}, 
the nominal trajectory $\alpha(t)$ and its negative version $-\alpha(t)$ are used
in place of ${\bar\alpha}^a(t)$ and ${\bar\alpha}^p(t)$.  
For a bouncing ball impacting without energy loss on a surface (located at position zero), 
the use of the mirror trajectory $-\alpha(t)$ can be justified within our framework 
observing that $-\alpha(t)$ corresponds to choosing ${\bar \alpha}^p$ before the nominal event time and
to ${\bar \alpha}^a$ after the nominal event time, so that an equivalent switching law to 
\eqref{eq:u_cl} is obtained. When the bouncing is not elastic, $-\alpha(t)$ is no longer
a good representative of the extended behavior and needs to be corrected. 
Indeed, in \cite{BiWoHeSaNi12C_TrackingDissipativeImpacts}, the case of non-elastic impact is considered
and the mirrored nominal trajectory is corrected via the use of the impact restitution coefficient. 
Again, this can be interpreted as the need to obtain (an estimate of) the extended nominal trajectories 
${\bar\alpha}^a(t)$ and ${\bar\alpha}^p(t)$ for properly defining the notion of tracking error 
to deal with the difference between the nominal and perturbed event times.
A similar remark applies for the mirroring technique presented in \cite{FoTeZa13J_FollowTheBouncingBall}.
\hfill$\Box$

Our goal in this section is to discuss
why \eqref{eq:u_cl} is a suitable choice to design a trajectory tracking controller 
(assuming $\xi$ is of infinite extent) and indicate how the time-varying gain $K$ can be designed.
Using \eqref{eq:u_cl}, we obtain the following closed-loop ante- and post-event vector fields 
\begin{align}
  f_{cl}^a(x,t) & := f^a (x, \mu + K(t) ( {\bar\alpha}^a(t) - x ), t) , 
  \label{eq:fa_cl} \\
  f_{cl}^p(x,t) & := f^p (x, \mu + K(t) ( {\bar\alpha}^p(t) - x ), t) .
  \label{eq:fp_cl}
\end{align}
By construction, the resulting hybrid system (with no inputs) has $\alpha$ as nominal trajectory and
consequently the nominal switching time remains $\te$.
The sensitivity analysis developed in Section~\ref{sec:sensitivity} leads to
the following state matrices for the ante- and post-event linearization of 
\eqref{eq:fa_cl}-\eqref{eq:fp_cl}:
\begin{align}
  A^a_{cl}(t) & = A^a(t) - B^a(t) K(t) , \\
  A^p_{cl}(t) & = A^p(t) - B^p(t) K(t) .
\end{align}
The input matrices $B^a_{cl}$ and $B^p_{cl}$
are zero as \eqref{eq:fa_cl} and \eqref{eq:fp_cl} have no input.
Finally, the gain $H$, computed using \eqref{eq:Hstep}, 
is by construction equal to the one associated to the nominal open-loop trajectory $\xi = (\alpha, \mu)$.
This concludes the derivation of the extended linearization as discussed in Proposition~\ref{prop:approximation}
for the closed-loop dynamics \eqref{eq:fa_cl}-\eqref{eq:fp_cl}.

In virtue of Proposition~\ref{prop:approximation} and, in particular, of the approximation \eqref{eq:xeps_alpha_plus_z}, 
we expect to be able to shape the {\em local} behavior of the closed-loop response of the hybrid system about the nominal trajectory $\xi$
by choosing suitably the matrix gain $K(t)$ in \eqref{eq:u_cl}. We clearly expect that this will be related to controllability-like assumptions 
on the ante- and post-event linearizations associated to the time-varying matrices $A^a$, $B^a$, $A^p$, $B^p$, as well as
the jump gain $H$ given in \eqref{eq:Hstep}. Here we limit ourselves to expose the main idea of using a modification
of the standard linear quadratic regulator (LQR) problem in order to compute the gain $K$ and leave to future investigation 
the task of filling in the gaps by providing a mathematical proof that the approach will in fact be effective
and demonstrate numerically the strategy on suitable examples.

Over a finite horizon, we consider the minimization of the quadratic cost functional
\begin{align}
   \frac{1}{2} \int_0^T z^T Q \,z + v^T R \,v  \, d s + \frac{1}{2} z(T)^T P_T z(T)
   \label{eq:JLQRcost}
\end{align}
\noindent subject to the jump linear dynamics 
\begin{align}
      z  & = z_0,                   &   t & = t_0,          \\                
 \dot z  & = A^a(t) z + B^a(t) v,   &   t & \in [t_0,\te],  \\
    z^+  & = z^- + H z^-,           &   t & = \te,          \\
 \dot z  & = A^p(t) z + B^p(t) v,   &   t & \in [\te,T].     
 \label{eq:JLQRdyn_eq4} 
\end{align}
with $Q$, $R$, and $P_T$ being (possibly time-varying) positive definite symmetric matrices.

The solution of the above optimal control
problem is given by the closed-loop feedback law 
\begin{align}
  v = - R^{-1} B^T P(t) z =: - K(t) z
  \label{eq:v_cl}
\end{align}
where $P(t)$ is the solution of the following Riccati differential equation with jumps
\begin{align}
        P   & =  P_T,                                       &   t & =  T,           \notag\\
   -\dot P  & = A^p P + P A^p - P B^p R^{-1} (B^p)^T P + Q, &   t & \in [\te, T],   \notag\\
       P^-  & = (I + H)^T P^+ (I + H),                      &   t & =   \te,        \notag\\
  -\dot P   & = A^a P + P A^a - P B^a R^{-1} (B^a)^T P + Q, &   t & \in [t_0, \te]. \notag 
\end{align}
The proof that \eqref{eq:v_cl} is indeed the solution 
to the optimal control problem \eqref{eq:JLQRcost}-\eqref{eq:JLQRdyn_eq4} 
is readily obtained by observing that \eqref{eq:JLQRcost}-\eqref{eq:JLQRdyn_eq4}
is equivalent to the minimization of 
\begin{align}
   \frac{1}{2} \int_0^\te z^T Q \, z + v^T R \, v  \, d s + \frac{1}{2} z(\te)^T P^- z(\te)
\end{align}
\noindent subject to the linear dynamics 
\begin{align}
       z & = z_0,                   &   t & = t_0,          \\                
 \dot z  & = A^a(t) z + B^a(t) v,   &   t & \in [t_0,\te],  
\end{align}
where $P^- = (I + H)^T P^+ (I + H)$ so that 
$1/2 (z^-)^T P^- z^- = 1/2 (z^+)^T P^+ z^+$ is the value function associated to the LQR problem
\begin{align}
   \frac{1}{2} \int_\te^T z^T Q \, z + v^T R \, v  \, d s + \frac{1}{2} z(T)^T P_T z(T)
\end{align}
\noindent subject to the linear dynamics 
\begin{align}
      z  & = z^+,                   &   t & = \te,          \\                
 \dot z  & = A^p(t) z + B^p(t) v,   &   t & \in [\te,T].
\end{align}
In case the nominal trajectory $\alpha$ becomes constant (resp., periodic) after the finite horizon time $T$,
$P_T$ can be initialized to the corresponding algebraic (resp., periodic) solution of the Riccati differential equation
to obtain a gain $K$ defined over the infinite horizon. 
It is left for future investigation the study of the effect on the optimal gain $K$ of the reset map $P^- = (I + H)^T P^+ (I + H)$,
found in the Riccati differential equation with jumps associated to the optimal control problem
\eqref{eq:JLQRcost}-\eqref{eq:JLQRdyn_eq4}. Our claim is that, as a consequence of this reset, the gain $K$
will be reduced before the impact event. 

\section{Conclusions}
\label{sec:conclusions}

This paper addresses the sensitivity analysis 
of hybrid systems with discontinuous state trajectories.
We developed a novel notion of error to obtain, at each instant of time, 
a first-order approximation of the change in a trajectory due to small
changes in initial conditions and inputs.
This notion of error naturally finds application in the local tracking problem 
of a time-varying reference trajectory of a hybrid system.
We outlined how the standard linear trajectory tracking control
for nonlinear systems --based on linear quadratic regulator (LQR) theory 
to compute the optimal feedback gain-- can be generalized for hybrid systems. 
We highlighted the connection between the switching linear feedback law that we 
propose with the idea of trajectory mirroring recently appeared in the literature.

The notion of error developed in this paper opens the possibility 
of further developing perturbation analysis in the context of hybrid systems. 
Our current efforts are directed toward the development of a second-order approximation
that will find application in the context of numerical optimal control
\cite{SaHaAg13J_OptCtrlLieGroups, HaSa06C_BarrierMethodOptTrajFunctionals, Ha02C_ProjectionOperator}.

\bibliographystyle{plain}
\bibliography{CDC2014_Sensitivity}

\begin{thebibliography}{10}

\bibitem{BiWoHeNi13J_TrackCtrlHSStateTriggeredJumps}
J.~J.~B. Biemond, N.~van~de Wouw, W.~P. M.~H. Heemels, and H.~Nijmeijer.
\newblock {Tracking Control for Hybrid Systems With State-Triggered Jumps}.
\newblock {\em IEEE Transactions on Automatic Control}, 58(4):876--890, April
  2013.

\bibitem{BiWoHeSaNi12C_TrackingDissipativeImpacts}
J.J.B. Biemond, N.~van~de Wouw, W.P.M.H. Heemels, R.G. Sanfelice, and
  H.~Nijmeijer.
\newblock {Tracking control of mechanical systems with a unilateral position
  constraint inducing dissipative impacts}.
\newblock In {\em Proc. of the 51st IEEE Conference on Decision and Control
  (CDC)}, pages 4223--4228. Ieee, December 2012.

\bibitem{BR99B_NonsmoothMechanics}
B.~Brogliato.
\newblock {\em {Nonsmooth Mechanics. Models, Dynamics and Control}}.
\newblock Springer Verlag, 2nd edition, 1999.

\bibitem{BeBuChKo08B_PiecewiseSmoothDynSystems}
M.~di~Bernardo, C.~Budd, A.R. Champneys, and P.~Kowalczyk.
\newblock {\em {Piecewise-smooth Dynamical Systems}}.
\newblock Springer, 2008.

\bibitem{FoTeZa13J_FollowTheBouncingBall}
F.~Forni, A.R Teel, and L.~Zaccarian.
\newblock {Follow the Bouncing Ball: Global Results on Tracking and State
  Estimation With Impacts}.
\newblock {\em IEEE Transactions on Automatic Control}, 58(6):1470--1485, 2013.

\bibitem{GaFeBa99J_SensitivityFunctionsHybridSystems}
S.~Gal\'{a}n, W.F. Feehery, and P.I. Barton.
\newblock {Parametric sensitivity functions for hybrid discrete/continuous
  systems}.
\newblock {\em Applied Numerical Mathematics}, 31(1):17--47, September 1999.

\bibitem{GaMePo08C_TrajTrackLinearHS}
S.~Galeani, L.~Menini, and A.~Potini.
\newblock {Trajectory tracking in linear hybrid systems: An internal model
  principle approach}.
\newblock In {\em 2008 American Control Conference}, pages 4627--4632. Ieee,
  June 2008.

\bibitem{GoSaTe12B_HybridSystems}
R~Goebel, R.G. Sanfelice, and A.R. Teel.
\newblock {\em {Hybrid Dynamical Systems}}.
\newblock Princeton University Press, 2012.

\bibitem{Ha02C_ProjectionOperator}
J.~Hauser.
\newblock {A Projection Operator Approach to the Optimization of Trajectory
  Functionals}.
\newblock In {\em Proceedings of the 15th IFAC World Congress}, Barcelona,
  Spain, 2002.

\bibitem{HaSa06C_BarrierMethodOptTrajFunctionals}
J.~Hauser and A.~Saccon.
\newblock {A Barrier Function Method for the Optimization of Trajectory
  Functionals with Constraints}.
\newblock In {\em Proceedings of the 45th IEEE Conference on Decision and
  Control}, pages 864--869. IEEE, 2006.

\bibitem{HiPa00J_SensitivityHybridSystems}
I.A. Hiskens and M.A. Pai.
\newblock {Trajectory Sensitivity Analysis of Hybrid Systems}.
\newblock {\em IEEE Transactions on Circuits and Systems - Part I: Fundamental
  Theory and Applications}, 47(2):204--220, 2000.

\bibitem{Kh02B_NLSystems}
H.K. Khalil.
\newblock {\em {Nonlinear Systems}}.
\newblock Prentice Hall, 3rd edition, 2002.

\bibitem{Kr85B_SeminalSensitivityAccordingToMombaur}
P.~Kr\"{a}mer-Eis.
\newblock {Ein Mehrzielverfahren zur numerischen Berechnung optimaler
  Feedback-Steuerungen bei beschr\"{a}nkten nichtlinearen Steuerungsproblemen}.
\newblock In {\em Bonner Mathematischen Schriften 166}. 1985.

\bibitem{LeinevandeWouw2008BOOK}
R.~I. Leine and N.~van~de Wouw.
\newblock {\em Stability and Convergence of Mechanical Systems with Unilateral
  Constraints}, volume~36 of {\em Lecture Notes in Applied and Computational
  Mechanics}.
\newblock Springer Verlag, Berlin, 2008.

\bibitem{LeCaVr00J_BifurcNLDiscontinuousSys}
R.I. Leine, D.H. van Campen, and B.L. van~de Vrande.
\newblock {Bifurcations in Nonlinear Discontinuous Systems}.
\newblock {\em Nonlinear Dynamics}, 23:105--164, 2000.

\bibitem{LeWo08J_UnifConvofMonotoneMDI}
R.I. Leine and N.~{van de Wouw}.
\newblock {Uniform Convergence of Monotone Measure Differential Inclusions:
  With Application To the Control of Mechanical Systems With Unilateral
  Constraints}.
\newblock {\em International Journal of Bifurcation and Chaos},
  18(05):1435--1457, May 2008.

\bibitem{MeTo01J_AsympTrackPeriodicTrajsNonsmothImpacts}
L.~Menini and A.~Tornamb\`{e}.
\newblock {Asymptotic Tracking of Periodic Trajectories for a Simple Mechanical
  System Subject to Nonsmooth Impacts}.
\newblock {\em IEEE Transactions on Automatic Control}, 46(7):1122--1126, 2001.

\bibitem{Mo01D_StabOptOpenLoopWalkRobots}
K.D. Mombaur.
\newblock {\em {Stability Optimization of Open-loop Controlled Walking
  Robots}}.
\newblock PhD thesis, Heidelberg University, Ruperto Carola, 2001.

\bibitem{MoGr09J_HybridInvManifoldSystemsWithImpulseEffects}
B.~Morris and J.W. Grizzle.
\newblock {Hybrid Invariant Manifolds in Systems With Impulse Effects With
  Application to Periodic Locomotion in Bipedal Robots}.
\newblock 54(8):1751--1764, 2009.

\bibitem{No91J_NonPeriodMotionGrazingImpactOscillator}
A.B. Nordmark.
\newblock {Non-Periodic Motion caused by Grazing Incidence in an Impact
  Oscillator}.
\newblock {\em Journal of Sounds and Vibration}, 145(2):279--297, 1991.

\bibitem{Pa01J_CtrlContactProblemConstrainedEL}
P.R. Pagilla.
\newblock {Control of Contact Problem in Constrained Euler–Lagrange Systems}.
\newblock {\em IEEE Transactions on Automatic Control}, 46(10):1507--1509,
  2001.

\bibitem{PaYu01C_ExpStudyPlanarImpactRobotManipulator}
P.R. Pagilla and B.~Yu.
\newblock {An Experimental Study of Planar Impact}.
\newblock In {\em Proc. of the IEEE International Conference on Robotics \&
  Automation (ICRA)}, pages 3943--3948, 2001.

\bibitem{RoLeSe06J_BounceJuggling}
R.~Ronsse, P.~Lef\`{e}vre, and R.~Sepulchre.
\newblock {Sensorless Stabilization of Bounce Juggling}.
\newblock {\em IEEE Transactions on Robotics}, 22(1):147--159, 2006.

\bibitem{Ro67J_GenSensitivityEqsDiscSystems}
E.N. Rozenvasser.
\newblock {General sensitivity equations of discontinuous systems}.
\newblock {\em Automation and Remote Control}, pages 400--404, 1967.

\bibitem{SaHaAg13J_OptCtrlLieGroups}
A.~Saccon, J.~Hauser, and A.~P. Aguiar.
\newblock {Optimal Control on Lie Groups: The Projection Operator Approach}.
\newblock {\em IEEE Transactions on Automatic Control}, 58(9):2230--2245,
  September 2013.

\bibitem{ScWiSc96C_ParEstimDiscDescriptorModel}
R.~von Schwerin, M.~Winckler, and V.~Schulz.
\newblock {Parameter estimation in discontinuous descriptor models}.
\newblock In {\em Proceedings of the IUTAM Symposium on Optimization of
  Mechanical Systems}, pages 269--276, 1996.

\bibitem{WeGrKo03J_HybridZeroDynPlanarBipedWalker}
E.R. Westervelt, J.W. Grizzle, and D.E. Koditschek.
\newblock {Hybrid Zero Dynamics of Planar Biped Walkers}.
\newblock {\em IEEE transactions on automatic control}, 48(1):42--56, 2003.

\end{thebibliography}

\end{document}